\documentclass{gtart_h}

\def\ifplaintex{\expandafter\ifx\csname documentclass\endcsname\relax}

\def\gtp{{\mathsurround=0pt\it $\cal G\mskip-2mu$eometry \&\ 
$\cal T\!\!$opology $\cal P\!$ublications}}  

\def\recd{{\small Received:\qua\receiveddate\ifx\reviseddate\relax
\else\qquad Revised:\qua\reviseddate\fi\par}} 

\def\lognumber#1{\def\thelognumber{#1}}
\def\volumenumber#1{\def\thevolumenumber{#1}}
\def\volumeyear#1{\def\thevolumeyear{#1}}
\def\papernumber#1{\def\thepapernumber{#1}}
\def\pagenumbers#1#2{\def\startpage{#1}\def\finishpage{#2}}
\def\published#1{\def\publishdate{#1}}

\def\received#1{\def\receiveddate{#1}}
\def\revised#1{\def\reviseddate{#1}}
\def\accepted#1{\def\accepteddate{#1}}

\def\asciiaddress#1{\def\theasciiaddress{#1}}

\long\def\asciiabstract#1{\long\def\theasciiabstract{#1}}

\let\\\par\let\thelognumber\relax\let\thevolumenumber\relax
\let\thepapernumber\relax\let\thevolumeyear\relax\let\startpage\relax
\let\finishpage\relax\let\publishdate\relax\let\receiveddate\relax
\let\reviseddate\relax\let\accepteddate\relax\let\theasciititle\relax
\let\theasciiauthors\relax\let\theasciiaddress\relax
\let\theasciiabstract\relax

\let\theasciiemail\relax

\ifplaintex
\font\logobig=cmssbx10 scaled 3836
\font\logomed=cmssbx10 scaled 2557
\else
\font\logobig=cmssbx10 scaled 4200
\font\logomed=cmssbx10 scaled 2800
\fi

\long\def\makeagttitle{   
\count0=\startpage
\agt\hfill      
\hbox to 45truept{\vbox to 0pt{\vglue -13truept{\logomed A\kern -.37em{\logobig 
T}\kern -.38em G}\vss}\hss}
\break
{\small Volume \thevolumenumber\ (\thevolumeyear)
\startpage--\finishpage\nl
Published: \publishdate}

\vglue .25truein

{\parskip=0pt\leftskip 0pt plus
1fil\def\\{\par\smallskip}{\Large\bf\thetitle}\par\medskip} \vglue
0.05truein

{\parskip=0pt\leftskip 0pt plus 1fil\def\\{\par}{\sc\theauthors}
\par\medskip}%
 
\vglue 0.03truein 

{\small\leftskip 25truept\rightskip 25truept{\bf Abstract}\stdspace\theabstract

{\bf AMS Classification}\stdspace\theprimaryclass
\ifx\thesecondaryclass\relax\else; \thesecondaryclass\fi\par
{\bf Keywords}\stdspace \thekeywords\par}\vglue 7truept

}   

\ifplaintex
\font\phead=cmsl9 scaled 950
\font\pnum=cmbx10 scaled 913
\font\pfoot=cmsl9 scaled 950
\headline{\vbox to 0pt{\vskip -4.5mm\line{\small\phead\ifnum
\count0=\startpage ISSN 1472-2739 (on-line) 1472-2747 (printed)
\hfill {\pnum\folio}\else\ifodd\count0\def\\{ }%
\ifx\theshorttitle\relax\thetitle\else\theshorttitle\fi\hfill{\pnum\folio}
\else\def\\{ and }{\pnum\folio}\hfill\ifx\theshortauthors\relax\theauthors
\else\theshortauthors\fi\fi\fi}\vss}}
\footline{\vbox to 0pt{\vglue 0mm\line{\small\pfoot\ifnum\count0=\startpage
\copyright\ \gtp\hfill\else
\agt, Volume \thevolumenumber\ (\thevolumeyear)\hfill\fi}\vss}}
\else
\font=cmsl9 scaled 1050
\font\lnum=cmbx10 
\font=cmsl9 scaled 1050
\makeatletter
\def\@oddhead{{\small\lhead\ifnum\count0=\startpage ISSN 1472-2739 
(on-line) 1472-2747 (printed)\hfill {\lnum\number\count0}\else\ifodd\count0
\def\\{ }\ifx\theshorttitle\relax \thetitle \else\theshorttitle\fi\hfill
{\lnum\number\count0}\else\def\\{ and }{\lnum\number\count0}
\hfill\ifx\theshortauthors\relax 
\theauthors\else\theshortauthors\fi\fi\fi}}\def\@evenhead{\@oddhead}
\def\@oddfoot{\small\lfoot\ifnum\count0=\startpage\copyright\ \gtp\hfill\else
\agt, Volume \thevolumenumber\ (\thevolumeyear)\hfill\fi}
\def\@evenfoot{\@oddfoot}
\makeatother
\fi
\let\maketitlepage\makeagttitle
\let\makeshorttitle\maketitlepage
\let\maketitle\maketitlepage

\newwrite\gtoutfile
\long\gdef\makeheadfile{  
{\def\\{, }\def\s{ }
\immediate\openout\gtoutfile head.xxx
\immediate\write\gtoutfile{Proxy-for: \ifx\theasciiauthors\relax
\theauthors\else\theasciiauthors\fi\s<\ifx\theasciiemail\relax\theemail\else\theasciiemail\fi>}
\immediate\write\gtoutfile{\noexpand\\}
\immediate\write\gtoutfile{Authors: \ifx\theasciiauthors\relax
\theauthors\else\theasciiauthors\fi}
{\def\\{ }\immediate\write\gtoutfile{Title: \ifx\theasciititle\relax
\thetitle\else\theasciititle\fi}}
\immediate\write\gtoutfile{Subj-class: GT or SG, GR etc}
\immediate\write\gtoutfile{MSC-class: \theprimaryclass\ifx\thesecondaryclass\relax\else, \thesecondaryclass\fi}
\immediate\write\gtoutfile{Journal-ref: Algebr. Geom. Topol. \thevolumenumber\s
(\thevolumeyear) \startpage-\finishpage}
\immediate\write\gtoutfile{Comments: Published by Algebraic and
Geometric Topology at}
\immediate\write\gtoutfile{\s\s\s  http://www.maths.warwick.ac.uk/agt/AGTVol\thevolumenumber/agt-\thevolumenumber-\thepapernumber.abs.html}
\immediate\write\gtoutfile{\noexpand\\}
\immediate\write\gtoutfile{}
\ifx\theasciiabstract\relax
\immediate\write\gtoutfile{\theabstract}\else
\immediate\write\gtoutfile{\theasciiabstract}\fi
\immediate\write\gtoutfile{}
\immediate\write\gtoutfile{\noexpand\\}
\immediate\write\gtoutfile{}
\immediate\closeout\gtoutfile}}  

\def\maketitlepage{\makeagttitle\makeheadfile}
\let\makeshorttitle\maketitlepage
\let\maketitle\maketitlepage

\lognumber{43}
\volumenumber{5}
\volumeyear{2005}
\papernumber{43}
\pagenumbers{1051}{1073}
\received{29 March 2005} 
\revised{30 June 2005}
\accepted{4 August 2005}
\published{29 August 2005}

\usepackage{amsmath,amssymb,color,epsfig,delarray}
\usepackage[all]{xy}

\def\figref#1{\hyperlink{#1anchor}{Figure~\ref*{#1}}}
\def\anchor#1{\noindent\hypertarget{#1anchor}{\smash{$\phantom{99}$}}\newline}

\def\circit#1#2#3{\vbox to 0pt{\kern-#2pt{\hbox to 
0pt{\kern#1pt{$_\circ$}\hss}\vss}}#3}
\def\Dcirc{\circit{3.5}{11.9}D}

\begin{document}
\newtheorem {lemma}{Lemma}
\newtheorem {defi}{Definition}
\newtheorem {prop}{Proposition}
\newtheorem {cor}{Corollary}
\newtheorem {theorema}{Theorem}
\newtheorem {claim}{Claim}
\newtheorem {fact}{Fact}
\newtheorem {remark}{Remark}
\newtheorem {recall}{Recall}
\newtheorem {example}{Example}

\title{Non-singular graph-manifolds of dimension 4}
\author{A. Mozgova}

\address{Laboratoire d'analyse non lin\'eaire et g\'eom\'etrie,
Universit\'e d'Avignon\\33, rue Louis Pasteur,
84000 Avignon, France{\rm \qua and}\\
Laboratoire Emile Picard, UMP 5580,
Universit\'{e} Paul Sabatier\\118, route de Narbonne, 31062 Toulouse,
France} 
\asciiaddress{Laboratoire d'analyse non lineaire et geometrie,
Universite d'Avignon\\33, rue Louis Pasteur,
84000 Avignon, France, and\\
Laboratoire Emile Picard, UMP 5580,
Universite Paul Sabatier\\118, route de Narbonne, 31062 Toulouse,
France}

\email{mozgova@univ-avignon.fr}

\begin{abstract}
A compact $4$-dimensional manifold is a non-singular graph-manifold if it 
can be obtained by the glueing $T^2$-bundles over compact surfaces 
(with boundary) of negative Euler characteristics. If none
of glueing diffeomorphisms respect the bundle structures, 
the graph-structure is called reduced. We prove that any homotopy equivalence 
of closed oriented $4$-manifolds with reduced nonsingular
graph-structures is homotopic to a diffeomorphism
preserving the structures. 
\end{abstract}

\asciiabstract{%
A compact 4-dimensional manifold is a non-singular graph-manifold if it 
can be obtained by the glueing T^2-bundles over compact surfaces 
(with boundary) of negative Euler characteristics. If none
of glueing diffeomorphisms respect the bundle structures, 
the graph-structure is called reduced. We prove that any homotopy equivalence 
of closed oriented 4-manifolds with reduced nonsingular
graph-structures is homotopic to a diffeomorphism
preserving the structures.} 



\maketitle

\section*{Introduction}
\addcontentsline{toc}{section}{Introduction}

In the paper \cite{Wald-gr}, Waldhausen introduced a
class of orientable $3$-manifolds called graph-manifolds which can
be obtained by glueing  blocks that are Seifert manifolds along
homeomorphisms of their boundary tori. These manifolds are not
always sufficiently large, but for them one can introduce a notion
of reduced graph-structure (i.e.\ a structure in which no family of
neighboring blocks can be replaced by a single block), and then, 
with a few explicit exceptions, the existence of a homeomorphism 
between two $3$-dimensional graph-manifolds with reduced
graph-structures implies the existence of a homeomorphism
respecting reduced graph-structures, which leads to a
classification of such $3$-manifolds.

$3$-dimensional graph-ma\-ni\-folds are important because 
they naturally arise as the boundary of resolved
isolated complex singularities of polynomial maps
$(\mathbb{C}^2, 0) \rightarrow (\mathbb{C}, 0)$ \cite{EisNeu},
as the surfaces of constant energy of integrable hamiltonian 
systems with two degree of freedom \cite{FomZi}, and
as $3$-manifolds admitting an injective $F$-structure (a
generalization of an injective torus action) \cite{Rong}.

Our goal is to study a class of smooth four-dimensional manifolds
generalizing three-di\-men\-si\-o\-nal graph-manifolds (with
blocks without singular fibers) and having fundamental groups of
exponential growth (hence, to which the high-dimensional techniques
do not apply).

\medskip

\textbf{Definition}\qua 1)\qua A (nonsingular) block is a
$T^2$-bundle over a compact surface (with boundary) of negative 
Euler characteristic. 

2)\qua A (nonsingular) graph-manifold structure on a manifold
is a decomposition as a union of blocks, glued by diffeomorphisms 
of the boundary.
\medskip

Note that the boundary of a block has the structure of a $T^2$-bundle
over a circle.

\medskip
\textbf{Definition}\qua A graph-manifold structure is
{\it reduced} if none of the glueing maps are isotopic 
to fiber-preserving maps of $T^2$-bundle.
\medskip

Any graph-structure gives rise to a reduced one by 
forming blocks glued by bundle maps into larger blocks. 

\medskip
\textbf{Main theorem}\qua {\sl Any homotopy equivalence 
of closed oriented $4$-manifolds with reduced nonsingular
graph-structures is homotopic to a diffeomorphism
preserving the structures.}
\medskip

The text is organized as follows. Section~1 
contains the main technical result.
A standard fact about two incompressible surfaces in an
orientable irreducible 3-manifold is that one can move one of them by 
isotopy in such a way that the new intersection becomes $\pi_1$-injective.
We provide a basis for doing a similar thing for $\pi_1$-injective
maps of $3$-manifolds into a $4$-manifold $W$ with
$\pi_2(W)=\pi_3(W)=0$ and for moving by (regular) homotopy.
The gain is the same: the intersection of images
of $3$-manifolds becomes completely visible in $\pi_1(W)$.
Section~2 contains a recapitulation of facts about
$T^2$-bundles over aspherical spaces. Section~3 
introduces four-dimensional non-singular graph-manifolds and proves the
main theorem. 


\section{3-dimensional $\pi_1$-injective submanifolds in 4-mani\-folds}

Manifolds here will be $C^{\infty}$. Denote the tangent map of $f$ by 
$df: TM \rightarrow TW$. An
{\it immersion} is a smooth map $f: M \rightarrow W$  such that at
every point of $M$ the derivative $df$ is an injective (linear)
map. The set of immersions is open in the
$C^{\infty}(M,W)$-topology (\cite{Munkres}, Theorem~3.10).
A {\it regular homotopy} is a homotopy through immersions. 
Any immersion $F: M \times I \rightarrow W$
such that $F \arrowvert_{M \times \{ 0 \} } =f_0$ and $F
\arrowvert_{M \times \{ 1 \} } =f_1$ gives a regular
homotopy between $f_0$ and $f_1$.
The converse is false: two embedded not
concentric circles in $\mathbb{R}^2$ are regularly homotopic, but their
embeddings can not be extended into an immersion $S^1 \times I \to
\mathbb{R}^2$.

\begin{figure}[ht]\anchor{reghomconstr}
\centering
\input{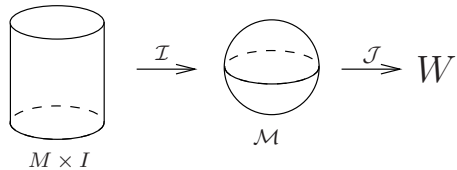}
\protect\caption{{\small Construction of a regular homotopy
}}\label{reghomconstr}
\end{figure}

\medskip{\bf Construction of regular homotopies}\qua
One particular method to construct a regular homotopy between two
immersions $f_0,\, f_1: M \rightarrow W$
is to immerse into $W$ not $M
\times I$, but the image of an isotopy. Precisely, suppose one has
a manifold $\mathcal{M}$ which contains the image of an isotopy
$\mathcal{I}$ between two embeddings $i_0: M \hookrightarrow
\mathcal{M},\ i_1: M \hookrightarrow \mathcal{M}$, i.e.\ there is a
map: $$ \mathcal{I}: M \times I \rightarrow \mathcal{M} \ \
\  \mbox{such that} \ \ \ \mathcal{I} \arrowvert_{M \times {0}}
=i_0, \ \mathcal{I} \arrowvert_{M \times {1}} =i_1$$ and
$\mathcal{I} \arrowvert_{M \times \{ t \} } \equiv \mathcal{I}_t $
is an embedding. Suppose also that there is an immersion
$\mathcal{J}: \mathcal{M} \rightarrow W$ such that $\mathcal{J I}
\arrowvert_{M \times {0}}= \mathcal{J} i_0 =f_0$ and $\mathcal{J
I} \arrowvert_{M \times {1}}= \mathcal{J} i_1 =f_1.$ Then the map
$\mathcal{J I}: M \times I \rightarrow W$ gives a smooth regular
homotopy between the immersions $f_0$ and $f_1$ (see
\figref{reghomconstr}).

\noindent We will refer to the above construction
by saying \textsl{``push $f_0(M)$ to $f_1(M)$ across $\mathcal{JI}
(M \times I)$''}.

\medskip{\bf Extension of immersions}\qua
Let $H: M \times I
\rightarrow W$ be a map such that $H \arrowvert_{M \times \{  0
\}}$ is an immersion and $\mathsf{dim} \, M + 1 < \mathsf{dim}\,
W$. The immersion $H \arrowvert_{M \times \{ 0 \}}: M \rightarrow
W$ determines a bundle injection $$A: T(M \times \{ 0 \})
\rightarrow (H^{\ast}TW) \arrowvert_{M \times \{ 0 \}}=(H
\arrowvert_{M \times \{ 0 \}})^{\ast}TW.$$
Suppose that this bundle injection can be extended to an
injection $A_I: T(M\times I) \rightarrow H^{\ast}TW$. Then, by
the Immersion Theorem (\cite{Smale}, \cite{Hirsch-immersion})
there exists an immersion $\mathcal{H}: M \times
I \rightarrow W$ which is homotopic to $H$, inducing the same
tangent bundle injection: $$\xymatrix{ T(M \times I) \ar[r]^{A_I}
\ar[rd]_{d \mathcal{H}} & H^{\ast}TW \ar@{=}[d] \\ &
\mathcal{H}^{\ast}TW}$$
The immersion $\mathcal{H}$ can be chosen in such a way that
$\mathcal{H} \arrowvert_{M \times \{ 0 \}}=H \arrowvert_{M \times
\{ 0 \}}$.

\medskip{\bf Main technical result}\qua
The following proposition is a generalization and a detailed
proof of the Proposition 2.B.2 of \cite{Stallings-book}, where very few
details of the proof are given.

\begin{prop}\label{reghom}
Let $W$ be a compact smooth oriented 4-dimensional manifold with
$\pi_2(W)=0$ and $M_1, \ M_2$ be compact oriented 3-manifolds with
$\pi_2(M_1)=\pi_2(M_2)=0$. Let $f_1:M_1 \rightarrow W$ be a
$\pi_1$-injective map and $f_2:M_2 \rightarrow W$ be a
$\pi_1$-injective embedding.

Then

\begin{itemize}
\item
$f_1$ is homotopic to a map $\widetilde{f}_1$ such that each
connected component of ${\widetilde{f}_1}^{-1}(f_2(M_2))$ is
$\pi_1$-injective in $M_1$;

\item if $\pi_3(W)=0$ and $M_1$ is irreducible,
then all the $S^2$-components of
${\widetilde{f}_1}^{-1}(f_2(M_2))$ can be eliminated by
homotopy of $f_1$;

\item in addition, if $f_1$ is an immersion, then the homotopies
can be made regular.

\end{itemize}
\end{prop}
\begin{proof}
Move $f_1$ by a small (regular if $f_1$ immersion) homotopy to make
it transverse to the submanifold $f_2(M_2)$.  
Then $F=f^{-1}_1(f_2(M_2))$ is a closed 2-dimensional surface which is 
embedded into $M_1$ (and immersed into $M_2$ if $f_1$ is immersion):
{\small
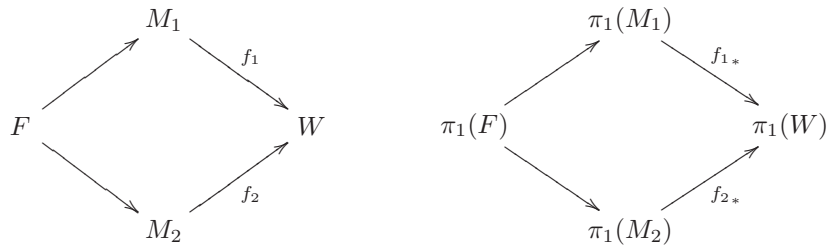
$$\xymatrix{ & M_1 \ar[rd]^{f_1} & & & \pi_1(M_1)
\ar[rd]^{{f_1}_{\ast}}  &
\\ \ \ \ \  F \ \ \ \  \ar[ru] \ar[rd]
& & \ \ \ \   W \ \ \ \   & \pi_1(F)
\ar[ru] \ar[rd] & & \pi_1(W) \\ &  M_2
\ar[ru]_{f_2} & & & \pi_1(M_2) \ar[ru]_{{f_2}_{\ast}} &   } $$}As
$f_2(M_2)$ is closed and $M_1$ is compact, the surface $F$ has only a 
finite number of connected components (\cite{Abraham-Robbin}, 
corollary~17.2(IV)).

\medskip


\medskip

\noindent {\bf Step~1\qua Construction of a map (resp.\ immersion) 
$\alpha: D^2 \times I \rightarrow W$ --- the image of the 
future homotopy (resp.\ regular homotopy)}

Suppose $ F \subset M_1$ is not $\pi_1$-injective, so $F$
is compressible in $M_1$, i.e.\ there exists an embedding 
$\beta : D^2 \rightarrow M_1$ such
that its boundary loop $\beta(\partial D^2)$ is not contractible in
$F$ but $\beta (D^2) \bigcap F = \beta(\partial D^2)$.

Consider the map $f_1 \beta: D^2 \rightarrow W$. As $F=f^{-1}(f_2(M_2))$
and $\beta(\partial D^2) \subset F$, hence $f_1 \beta (\partial D^2)
\subset f_2(M_2)$. Thus the map $f_1 \beta$ is in fact $$f_1 \beta:
(D^2, \partial D^2) \rightarrow (W, f_2(M_2)).$$ Note now that $\pi_2(W,
f_2(M_2))=0$ because for the embedding $f_2$ the corresponding map
induced in the fundamental groups ${f_2}_{\ast}: \pi_1(M_2)
\rightarrow \pi_1(W)$ is injective and the homotopy sequence of
the pair $(W, f_2(M_2))$ $$\xymatrix{\cdots \ar[r] & \pi_2(W)
\ar@{=}[d] \ar[r] & \pi_2(W, f_2(M_2)) \ar[r] & \pi_1(M_2)
\ar@{=}[d] \ar[r]^{{f_2}_{\ast}} & \pi_1(W) \ar[r] & \dots \\ & 0
&  & \pi_1(f_2(M_2)) && }$$ \noindent is exact. This implies that
$f_1 \beta$ is homotopic to a map $D^2 \to f_2(M_2)$ which can
be written as $f_2 \beta_2: D^2 \rightarrow f_2(M_2)$, i.e.\ there
exists a map $H: D^2 \times I \rightarrow W$ such that $H
\arrowvert_{D^2 \times \{ 0 \} }= f_1 \beta, \ H \arrowvert_{D^2
\times \{ 1 \} }= f_2 \beta_2$. More, the homotopy can be made in such a
way, that $\forall t \ \, H \arrowvert_{\partial D^2 \times \{ t \}}= 
H \arrowvert_{\partial D^2 \times \{ 0 \}}=f_1 \beta 
\arrowvert_{\partial D^2}$.

\medskip

{\bf Step~1.1\qua The case of ordinary homotopy}\qua In the
case when $f_1$ is just a {\it map} and we are interested in an 
{\it ordinary homotopy}, put $\alpha: = H$. As $D^4=D^3
\times I$ retracts on $D^3 \times \{ \frac{1}{2} \}=(D^2 \times I )
\times \frac{1}{2} $, we can say
that $\alpha$ extends to a map $D^4 \to W^4$: {\footnotesize
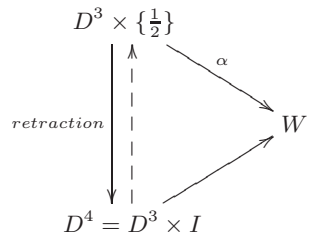
$$\xymatrix{D^3 \times \{ \frac{1}{2} \} \ \ \ar[rd]^{\alpha}
\ar@<-2ex>[dd]_{retraction} &
\\ & W
\\ D^4=D^3 \times I  \ar@{-->}[uu] \ar[ur] & }$$}
and move to the next step.

\medskip

{\bf Step~1.2\qua The case of regular homotopy}\qua In the
case when $f_1$ is an {\it immersion} (with trivial normal bundle since
everything is orientable) and we are looking for a
{\it regular homotopy}, let us show that this map $H$ can be
changed to an immersion.

As $H \arrowvert_{D^2 \times \{ 0 \} }= f_1 \beta: D^2
\rightarrow W$ is an immersion, its derivative $dH \arrowvert_{D^2
\times \{ 0 \} }= d( f_1 \beta) $ is correctly defined and gives
a bundle injection $$T(D^2 \times \{ 0 \}) \longrightarrow \bigl(
H \arrowvert_{D^2 \times \{ 0 \} } \bigr)^{\ast}TW= H^{\ast}TW
\arrowvert_{D^2 \times \{ 0 \}}.$$
Since $D^2 \times I$ retracts to $D^2 \times \{ 0 \}$,
this bundle injection extends to a bundle injection
$T(D^2 \times I) \longrightarrow H^{\ast}TW$ that on $D^2 \times \{ 1\}$
restricts to a subbundle of $(H \arrowvert_{D^2 \times I})^{\ast} T_{M_2}$.
Applying the Immersion Theorem gives an immersion
$\alpha: D^2 \times I \to W$ such that $\alpha
\arrowvert_{D^2 \times \{ 0 \} }= H \arrowvert_{D^2 \times \{ 0 \}
}= f_1 \beta$ and $\alpha(D^2 \times \{ 1 \} ) \subset
f_2(M_2)$.
Note that the immersion $\alpha$ is flat ($D^2 \times I$ and $W$
being oriented): there exists a map $$\alpha_{\varepsilon}:
D^4=(D^2 \times I) \times I \to W \ \ \ \mbox{such that} \ \ \
\alpha_{\varepsilon} \arrowvert_{(D^2 \times I) \times \{
\frac{1}{2} \}} \equiv \alpha.$$ As $\beta$ is flat (being an 
embedding), we have in $M_1$ an embedded $3$-disk
${\beta}_{\varepsilon}(D^2 \times \{ 0 \} \times I)$ which is the normal bundle
of $\beta(D^2)$. As $\alpha_{\varepsilon} \arrowvert_{(D^2
\times I) \times \{ \frac{1}{2} \}} \equiv \alpha$ and $\alpha
\arrowvert_{D^2 \times \{ 0 \} \times \frac{1}{2}}$ we can write
$\alpha_{\varepsilon} \arrowvert_{(D^2 \times \{ 0 \}) \times I} =
f_1 {\beta}_{\varepsilon}$.

\medskip

We will use the notation $\partial D^3= S^2_+ \cup S^2_-$ with $S^2_+=D^2
\times \{ 0 \}$ and $S^2_-=(D^2 \times \{ 1 \}) \cup (\partial
D^2 \times I)$.

\medskip

{\bf Step~2\qua Homotopy description}

\noindent As the result of the above construction we have:
\begin{itemize}
\item
an embedding of $2$-disk $\beta: D^2 \rightarrow M_1$ such that
$\beta(D^2) \cap F= \beta (\partial D^2)$,
\item
a map (respectively, a flat immersion) of $3$-disk $\alpha:
D^3=D^2 \times I \rightarrow W$ such that $\alpha
\arrowvert_{D^2 \times \{ 0 \}}=f_1 \beta$, $\alpha
\arrowvert_{D^2 \times \{ 1 \}} \subset f_2(M_2)$ and
$\alpha(\partial(D^2 \times \{ 0 \}))) =
\alpha(\partial(D^2 \times \{ 1 \})))$.
\end{itemize}

\begin{figure}[ht]
\centering
\input{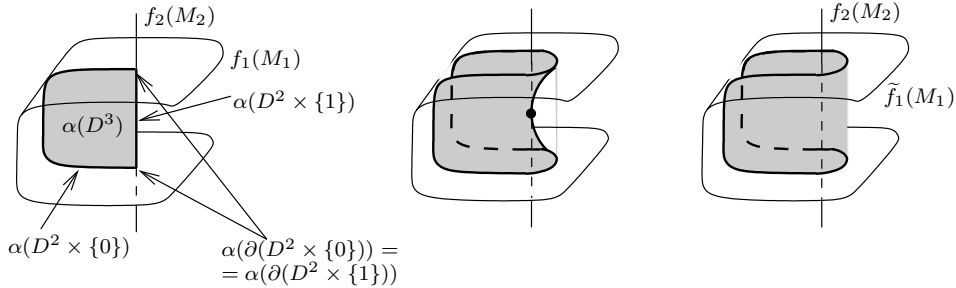}
\protect\caption{{\small Pushing ${f_1}
\arrowvert_{{\beta}_{\varepsilon}(D^2 \times I)}$ across
$\alpha_{\varepsilon}(D^3 \times I)$}}
\end{figure}

We will now change the map $f_1$
firstly by pushing $f_1 \arrowvert_{{\beta}_{\varepsilon}(D^2
\times I)}$ across $\alpha_{\varepsilon}(D^3 \times I)$ to a map
(respectively, an immersion) into $W$ whose image lies in
$f_2(M_2)$; secondly we compose it with the pushing along the
normal bundle of $f_2(M_2)$ in $W$ in such a way that:

\begin{itemize}

\item
in a small neighborhood $U \supset {\beta}_{\varepsilon}(D^2
\times I)$ the map (respectively, the immersion) $f_1$ is changed
by homotopy (resp.\ regular homotopy) to a map (respectively, an
immersion) $\widetilde{f}_1$ such that
${\widetilde{f}_1}^{-1}(f_2(M_2)):=F'$ is $F$ surgered on the disk
$\beta(D^2)$

\item  and the map $f_1$ does not change on the complement of $U$ in $M_1$.
\end{itemize}

\eject
\begin{figure}[h!]
\centering
\input{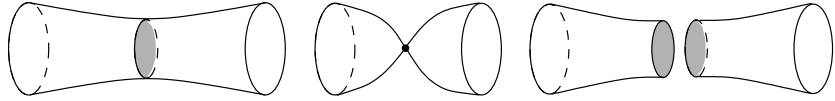}
\protect\caption{{\small Trace that the pushing of $f_1$ makes on
the intersection of images of $f_1$ and $f_2$}}
\end{figure}

The homotopy (resp.\ regular homotopy) works as follows.

\medskip

\noindent {\bf Step~3\qua Homotopy on a disk}

\noindent Let us decompose $\partial D^4$ as union of two $3$-discs
$S^3_+$ and $S^3_-$ with $S^3_+ \cap S^3_-=S^2$. Let $\mathcal{I}$ 
be an isotopy $D^3 \times I \to D^4$
that sends $3$-disk $S^3_+$ on $3$-disk $S^3_-$ as shown on
\figref{cub}.

\begin{figure}[ht]\anchor{cub}
\centering
\input{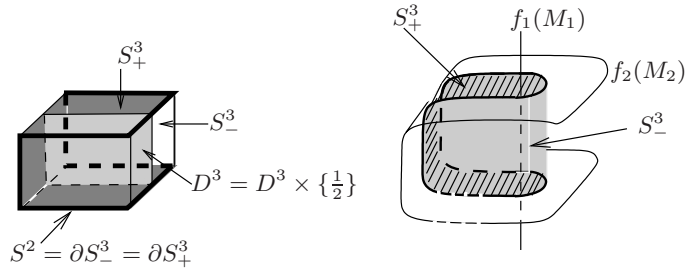}
\protect\caption{{\small Homotopy sending $S^3_+$ on $S^3_-$}}
\label{cub}
\end{figure}

\noindent Take now the composition $\alpha_{\varepsilon}
\mathcal{I}: D^3 \times I \to D^4$ ($\alpha_{\varepsilon}$ being the flat
extension of $\alpha$): it provides a homotopy
(respectively, a regular homotopy)
sending $\alpha_{\varepsilon} \arrowvert_{(D^2 \times \{ 0 \})
\times I} \equiv f_1 {\beta}_{\varepsilon}(D^3)$ into $f_2(M_2)$:
we ``push $f_1 \arrowvert_{\beta(D^2)}$ across
$\alpha_{\varepsilon} \mathcal{I} (D^4)$''.
Take then the composition of $\alpha_{\varepsilon} \mathcal{I}$
with the pushing out along the normal bundle of $f_2(M_2)$ in
$W^4$ (which is trivial, because $f_2(M_2)$ and $W^4$ are
orientable). At this moment the map $f_1$ will be changed not only
on the $3$-disc ${\beta}_{\varepsilon}(D^2 \times I)$, but on its small
neighborhood $U \subset M_1$.

\medskip

{\bf Step~4\qua Change $\alpha$ to make $\alpha(\Dcirc^3)$
miss $f_2(M_2)$}

{\bf Motivation}\qua If we want the homotopy
described on the previous step to create no new intersections of
$f_1(M_1)$ and $f_2(M_2)$, we have to make the image of the
interior of the disk $\alpha(\Dcirc^3)$ disjoint from $f_2(M_2)$.
We have $\alpha^{-1}(f_2(M_2))=S^2_- \cup G$, where $G$ are 
some closed surfaces.

Denote by $\Delta$ the union of $G$ and all components of $D^3 \setminus G$
that do not contain $\partial D^3$. Note some components of $G$ may be in
the interior of $\Delta$. Let $\hat G=\partial \Delta$. Since
$\Delta$ is an open subspace of a manifold, it is a manifold. 
Let us show that $\Delta$
is aspherical. If we show that $\pi_2(\Delta)=0$, it will
give us the asphericity: take the universal covering
$\widetilde {\Delta}$, we have $H_i(\widetilde{\Delta})=0,\ i \ge
3$ because it's an open $3$-manifold; then, by Whitehead's
theorem, $\pi_i(\widetilde{\Delta})=0,\ i \ge 3$, and we conclude
that $\pi_i(\Delta)=\pi_i(\widetilde{\Delta})=0$. Suppose
that $\pi_2(\Delta) \ne 0$. Then, by the Sphere Theorem, there
exists an embedded $S^2 \hookrightarrow \Delta$ representing a
non-trivial element in $\pi_2(\Delta)$. 
This $S^2$ bounds a ball in $D^3$. This ball must be contained in
$\Delta$, therefore $\pi_2(\Delta)=0$.

Note that $\Delta$ can be rather complex, for example, be a knot
complement.

\begin{figure}[ht]\anchor{preimage}
\centering
\input{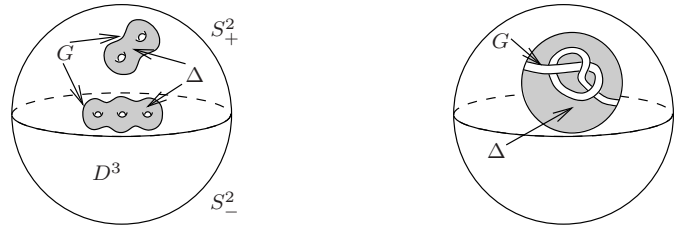}
\protect\caption{{\small Pre-image of $f_2(M_2)$ by $\alpha$: the
cases of $\Delta$ being handlebodies and a knot complement}}
\label{preimage}
\end{figure}

{\bf Extension of $\alpha \arrowvert_{\hat G}$ on
$\Delta$ and change $\alpha$}\qua Now, let's show that we can always
extend the map $\alpha \arrowvert_{\hat G}: \hat G \rightarrow f_2(M_2)$ to
a map $\alpha': \overline{\Delta} \rightarrow f_2(M_2)$. As $\hat G$
and $\overline{\Delta}$ are aspherical, it will be enough to
extend this map on the fundamental group of each component of $\Delta$. 
As $f_2$ is a
$\pi_1$-injective embedding, we have $\pi_1(M_2) \cong
\pi_1(f_2(M_2))$. $$\xymatrix{ \hat G \ar[r]^{\subseteq} \ar[d]_{\alpha
\arrowvert_{\hat G}} & \overline{\Delta} \ar[r]^{\subseteq}
\ar[d]_{\alpha \arrowvert_{\overline{\Delta} }}
\ar@{-->}[dl]_{\alpha'} & D^3 \ar[dl]^{\alpha} & & \pi_1(\hat G)
\ar[r]^{a} \ar[d]_{b} & \pi_1(\overline{\Delta}) \ar[r]^{d_1}
\ar[d]_{d} \ar@{-->}[dl]_{\alpha'_{\ast}} & 1 \ar[dl]^{d_2}
\\
f_2(M_2) \ar[r]_{\subseteq} & W &  & 1 \ar[r] & \pi_1(M_2)
\ar[r]_{c} & \pi_1(W) &} $$ \noindent We have $cb=d_2 d_1 a$, so
that $Im \ cb=1$. As $c$ is a monomorphism, it follows that $Im\
b=1$, thus, $b \equiv 1$. We can define the homomorphism
$\alpha'_{\ast}: \pi_1(\overline{\Delta}) \rightarrow \pi_1(M_2)
\cong \pi_1(f_2(M_2))$ as being the constant $1$, too.

So we can define a new map $\widetilde{\alpha}: D^3 \longrightarrow W$
as follows:
\[ \widetilde{\alpha}=
\begin{array}
\{{rl}. \alpha & \mbox{on } D^3 \backslash \Delta
\\
\alpha' & \mbox{on } \Delta \end{array} \]
\noindent Now the whole image $\widetilde{\alpha}(\Delta)$ lies in
$f_2(M_2)$, hence we can push $\widetilde{\alpha}(D^3)$ off
$f_2(M_2)$ across $\alpha'(\Delta) \subset f_2(M_2)$ using the
normal bundle of $f_2(M_2)$ in $W$. We have a new map (which for
simplicity we still note by $\alpha$) $\alpha: D^3
\longrightarrow W$ such that the image of the interiour 
$\alpha(\Dcirc^3)$ is disjoint from $f_2(M_2)$.

\medskip

{\bf Step~5\qua The number of disks in $M_1$, on which the
homotopy of $f_1$ must be done, is finite} 

Suppose we made the homotopy of the map $f_1$ on one disk. Suppose
that the obtained surface $F$ (whose topological type has changed) is
still not $\pi_1$-injective in $M_1$. After the surgery the new
surface is still oriented, hence again there is an embedded
compressing disk in $M_1$, on which again one can do the surgery by
homotopy of $f_1$ etc. After each surgery the topological type of the
surface $F$ changes as follows: either the genus of one component of
$F$ decrements, or one component splits into two components, the sum
of genera of which is not greater than the genus of the original
component.
As $F$ is compact, the genera of all its components are finite,
and as it was pointed out before Step~1, the number of components
of $F$ is finite, hence, {\it the process will terminate after a
finite number of steps}. This is the advantage that we get from
replacing the homotopic information ($Ker \, g_{\ast} \ne 0$) by
the geometric information (there exists an embedded loop which is
trivialized by an embedded disk): the infinite kernel is killed in
a finite number of steps.

As the result we obtain a surface that is $\pi_1$-injective in
$M_1$, but which could contain spheres among its components.

\medskip

{\bf Step~6\qua Elimination of $S^2$-components 
provided $\pi_3(W)=0$ and $M_1$ is
irreducible} 

If the obtained surface $F$ contains $S^2$-components, then, as $M_1$
is irreducible and $\pi_2(M_2)=0$, every such $S^2$-component bounds
an embedded $3$-disk in $M_1$ and a homotopy $3$-disk in $M_2$: there
exist an embedding $\gamma_1 : D^3 \to M_1$ and a map $\gamma_2: D^3
\to M_2$ such that $f_1 \gamma_1(\partial D^3)=f_2 \gamma_2(\partial
D^3).$ Denote $f_1 \gamma_1(D^3)= S^3_+$ and $f_2
\gamma_2(D^3)=S^3_-.$

\begin{figure}[ht]\anchor{kills2}
\centering
\input{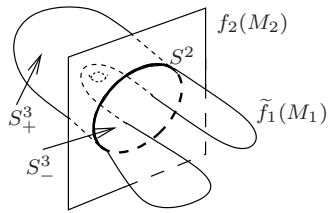}
\protect\caption{{\small Elimination of
$S^2$-components}}\label{kills2}
\end{figure}

\noindent If $\pi_3(W)=0$, then the map of $3$-sphere, whose image is $S^3_+
\bigcup S^3_-$, bounds a homotopy $4$-disk: there exists $\lambda:
D^3 \times I \to W$ such that $\lambda(D^2 \times \{ 0 \})=S^3_+,
\, \lambda(D^2 \times \{ 1\} = S^3_-$ and $\lambda(\partial D^3 \times I)=
S^3_+ \cap S^3_-$, and which, in addition, can be made an immersion on each
$D^3 \times \{ t \}$. In order to eliminate a chosen $S^2$-component of
$F$, push $\widetilde{f}_1 (\gamma_1(D^3))$ along $\lambda(D^3 \times I)
\subset W$ to make $S^3_+$ coincide with $S^3_-$ (see
\figref{kills2}), then, push it off $f_2(M_2)$ along the
normal bundle of $f_2(M_2)$ in $W$, then glue with
the map on $M_1 \setminus S^3_+$.
\end{proof}


\section{Torus bundles}

A {\it torus bundle} here will be a fiber bundle
$f: M \to B$ with fibers diffeomorphic to $T^2$, smooth if the base is a
smooth manifold.  The {\it monodromy} is
the action of $\pi_1(B)$ on $H_1$ of the fiber:
$$ \pi_1(B,b) \to Aut \, (H_1(f^{-1}(b)); \mathbb{Z}).$$
Choosing an identification of the fiber with $T^2$ (equivalently,
a basis for $H_1(T^2)$) identifies the automorphism group as $GL(2,
\mathbb{Z})$. The classifying map for a torus bundle is $B \to
B_{Diff(T^2)}$. There is a $2$-stage Postnikov decomposition $$K(\mathbb{Z}
\oplus \mathbb{Z}, 2) \to B_{Diff(T^2)} \to B_{GL(2, \mathbb{Z})}$$
(\cite{Hil}, ch.4, p.51). If $B$ is a surface
with non-empty boundary, this implies that a bundle is determined up to
isomorphism by the conjugacy class of its monodromy. However the
nontrivial $\pi_2$ in the classifying space shows bundles on surfaces are
not defined {\it rel boundary} by the monodromy.
A {\it fiber map} is a pair of maps,
one on the total spaces and one on the bases so that the diagram
$$\xymatrix{E_1 \ar@{-->}[r] \ar[d] & E_2 \ar[d] \\
B_1 \ar@{-->}[r] & B_2 }$$
\noindent commutes. 
{\it Bundle map} is a fiber map so that on coordinate charts 
it is given by function into the structure group.
As the inclusion $Diff(T^2) \hookrightarrow G(T^2)$ 
(the monoid of self-homotopy equivalences of torus) is homotopy equivalence
\cite{EE}, the existence of a bundle map between
$T^2$-bundles is equivalent to the existence of a fiber map 
inducing homotopy equivalence on the fibers.

A {\it fiber covering map} of bundles here will be a fiber map,
which is finite covering on fibers. The degree of the covering on
different fibers is clearly the same. 
If $B$ is aspherical, there exists a fiber covering map of
$T^2$-bundles with monodromies $\varphi_1, \, \varphi_2$ 
if and only if there exists a monomorphism 
$\alpha: \mathbb{Z} \oplus  \mathbb{Z} \to \mathbb{Z} \oplus \mathbb{Z}$
with $\alpha \, \varphi_1(\gamma) = \varphi_2(\gamma) \, \alpha$ for all
$\gamma \in \pi_1(B)$.  

\begin{lemma}\label{q1}
Let $f: E_1 \to E_2$ be $\pi_1$-injective map of $T^2$-bundles 
over aspherical spaces. Then $f$ is homotopic to a fiber covering 
map if and only if the induced map on $\pi_1$ sends the fiber subgroup 
of $\pi_1(E_1)$ into the fiber subgroup of $\pi_1(E_2)$. 
\end{lemma}
\begin{proof}
The $\pi_1$ condition is well-defined (independent of basepoints) because
the fiber defines a normal subgroup of $\pi_1$. A fiber covering map
clearly verifies the condition.

Now suppose $f$ sends the subgroup of fiber of $E_1$ into the
subgroup of fiber of $E_2$. Then it induces
an homomorphism on quotient groups. This gives a $\pi_1$-injective map of
the base spaces and a commutative up to homotopy diagram 
$$\xymatrix{ E_1 \ar^{f}[r]
\ar[d] & E_2 \ar[d] \\ G_1 \ar^{g}[r] & G_2 }$$
\noindent Let $E_2^{\ast}$ be the pullback of $E_2$ to $G_1$. Then $f$
factors as $f^{\ast}: E_1 \to E_2^{\ast}$ and 
a map $E_2^{\ast} \to E_2$ which is isomorphic on fibers, 
so it is sufficient to show that $f^{\ast}$ is
homotopic to a fiber covering map, which follows from the commutative
diagram of short exact sequences. 
\end{proof}

\begin{prop}\label{q2}
Suppose $E$ is homotopy equivalent to a $T^2$-fibration over a graph $G$.
Then this structure is unique up to homotopy unless $G \cong S^1$ and the
monodromy is conjugate to
$\bigl(
\begin{smallmatrix}
1 &  n \\ 0 & m
\end{smallmatrix}
\bigr)$.
\end{prop}
\begin{proof}
According to the lemma \ref{q1} it is sufficient to show that there is a
unique normal subgroup isomorphic to $\mathbb{Z} \oplus \mathbb{Z}$ and
with free quotient except when $G \cong S^1$ and the monodromy has the
specified form.

{\bf Case~0}\qua $G$ contractible, so $E \cong T^2$, and $\pi_1(E)
\cong \mathbb{Z} \oplus \mathbb{Z}$, then $\pi_1(E)$ is the only such
subgroup.

{\bf Case~1}\qua $G \cong S^1$. Then the fundamental group of the
fiber is the unique normal $( \mathbb{Z} \oplus \mathbb{Z})$-subgroup 
of $\pi_1(E)$ unless the monodromy has
an eigenvector with eigenvalue $1$. This shows the monodromy is conjugate
to $\bigl(
\begin{smallmatrix}
1 &  n \\ 0 & m
\end{smallmatrix}
\bigr)$. In this case either $\pi_1(E) \cong \mathbb{Z}^3$ or the
commutator subgroup is $\mathbb{Z}$, with quotient $\mathbb{Z}^2$. This
means $E$ is homotopic to an $S^1$-bundle over $T^2$. Taking
non-homotopic fibering $T^2 \to S^1$ induces non-homotopic
$T^2$-bundle structure on $E$.

{\bf Case~2}\qua $\pi_1(G)$ is non-abelian free group.
Consider the exact homotopy sequences of two $T^2$-bundles on $E$.
In the sequence of first bundle $$\xymatrix{ 0 \ar[r] &
\mathbb{Z} \oplus \mathbb{Z} \ar[r]^{\alpha} & \pi_1(E)
\ar[r]^{\beta}  & \pi_1(G) \ar[r] & 0 \\ 0 \ar[r] & \mathbb{Z} \oplus
\mathbb{Z} \ar[ur]_{j} & & & } $$
\noindent we have $Im \, j \subset Im\, \alpha=Ker \, \beta$, because
$\beta (Im \, j) \subset \pi_1(G)$ is an abelian normal subgroup, hence
trivial (theorem~2.10 of \cite{MagnKarrSoli}). Similarly, from the 
sequence of the second bundle we have $Im \, \alpha \subset Im \, j$,
hence these subgroups coincide in $\pi_1(E)$.
\end{proof}

\begin{cor}\label{coruni}
If a $4$-manifold is a $T^2$-bundle over a surface with boundary 
different from annulus and M\"{o}bius band, then this structure 
is unique up to bundle homotopy.
\end{cor}
\begin{proof}
A surface with boundary has the homotopy type of a graph.
\end{proof}

\begin{prop}\label{fiber-cover}
Let $f: E_1 \to E_2$ be a $\pi_1$-injective map between $T^2$-bundles
over aspherical surfaces. Then $f$ is homotopic to a fiber covering map
unless $E_1$ either comes from $S^1$-bundle over $3$-dimensional
$S^1$-bundle over aspherical surface (whose $\pi_1$ contains a normal 
$\mathbb{Z}$) or is $T^4$ or $T^2 \times K^2$.
\end{prop}
\begin{proof}
Denote the projection $p: E_2 \to B$ and $G:=Im \, (p_{\ast}f_{\ast})$.
Let $B_{2,G}$ be a covering of $B_2$ corresponding to $G$, 
$p_G: E_{2,G} \to B_{2,G}$ be the pullback of $p: E_2 \to B_2$ by 
$B_{2,G} \to B_2$, and $\hat{f} : E_1 \to E_{2,G}$ be a map covering $f$. 
$$\xymatrix{& E_1 \ar[dl]_{\hat{f}} \ar@{-->}[ddl] \ar[d]^{f} \\
E_{2,G} \ar[r]   \ar[d]_{p_G} & E_2 \ar[d]^{p}  \\ B_{2,G} \ar[r] & B_2  }$$
\noindent Denote the kernel of $\pi_1(E_1) \to \pi_1(B_1)$ 
by $K_1$ and the kernel of 
$(p_G \hat{f})_{\ast}: \pi_1(E_1) \to \pi_1(B_{2,G})$ by $K_2$. 
$$\xymatrix{K_1 \ar[r] & \pi_1(E_1) \ar[r] 
\ar[rd]_{(p_G \hat{f})_{\ast}} & \pi_1(B_1) \\
K_2 \ar[ur] & & \pi_1(B_{2,G})  }  $$
\noindent Both kernels are isomorphic to $\mathbb{Z} \oplus \mathbb{Z}$,
and lemma \ref{q1} implies that $f$ is homotopic to a fiber covering 
map if and only if $K_1=K_2$. 

{\bf Case~1}\qua $K_1 \cap K_2 \equiv \mathbb{Z} \oplus \mathbb{Z}$.

As $\pi_1$ of aspherical surface has no torsion, it means
$K_1=K_2$. It gives a map $B_1 \to B_{2,G}$ such that up to homotopy 
all the squares of the diagram commute
$$\xymatrix{\pi_1(E_1) \ar[rd]^{(p_G \hat{f})_{\ast}} 
\ar@/^2pc/[rr]^{f_{\ast}}
\ar[d] \ar[r]^{\hat{f}_{\ast}} & \pi_1(E_{2,G}) \ar[r]
\ar[d]^{(p_G)_{\ast}} & \pi_1(E_2) \ar[d]^{p_{\ast}} \\ 
\pi_1(B_1) \ar@{-->}[r] & \pi_1(B_{2,G}) \ar[r] &
\pi_1(B_2)   }$$ 
\noindent and hence a map $B_1 \to B_2$ which with $f$ gives a 
fiber covering map.

{\bf Case~2}\qua $K_1 \cap K_2 \equiv \mathbb{Z}$.

In this case $\pi_1(B_1)$ contains $\mathbb{Z}$ as normal
subgroup, hence $B_1$ is $T^2$, Klein bottle $K^2$, $S^1 \times I$ or
M\"{o}bius band.
As $(K_1 \cap K_2) \, \lhd \, \pi_1(E_1)$ and the monodromy acts by conjugation,
in this case the monodromy of $E_1 \to B_1$ preserves a curve in the
fiber. This curve is embeded because there is no torsion in $\pi_1(B_1)$
(and hence $K_1 / (K_1 \cap K_2) \cong \mathbb{Z}$). So that, $E_1$ is a
$S^1$-bundle over $3$-dimensional manifold $W$
$$\xymatrix{K_1 \cap K_2 \ar[r] & \pi_1(E_1) \ar[r] & \pi_1(W)}$$
\noindent and $W$ itself is $S^1$-bundle over $B_1$. Denote $K_3 \, \lhd
\, \pi_1(W)$ the corresponding fiber subgroup. The subgroup $K_2 / (K_1
\cap K_2) \cong \mathbb{Z}$ is normal in $\pi_1(W)$ with quotient
isomorphic to $\pi_1(B)$. The subgroups $K_3$ and $K_2 / (K_1 \cap K_2)$
coinside if and only if $f$ is fiber covering. 

If $K_3 \ne K_2 / (K_1 \cap K_2)$, refiber $W$ by $S^1$ with fiber
subgroup $K_2 / (K_1 \cap K_2)$. Together with $E_1 \to W$ it will
give another $T^2$-fibration of $E_1$, in which $f$ will be
fiber-covering. 

{\bf Case~3}\qua $K_1 \cap K_2 \equiv 1$.

In this case $\pi_1(B_1)$ contains a normal $\mathbb{Z} \oplus
\mathbb{Z}$, hence $B_1$ is $T^2$ or a Klein bottle, and $\pi_1(E_1)$ injects
into $\pi_1(B_1) \times \pi_1(B_{2,G})$. The monodromy of
$E_1 \to B_1$ is trivial, because in the diagram 
$$\xymatrix{0 \ar[d] & 0 \ar[d] & \\ K_1 \ar[r] \ar[d] & \pi_1(E_1)
\ar[r] \ar[d] & \pi_1(B_1) \ar@{=}[d] \\ \pi_1(B_{2,G}) \ar[r] &
\pi_1(B_1) \times \pi_1(B_{2,G}) \ar[r] & \pi_1(B_1)  }$$
\noindent the monodromy of the second line is trivial, the morphisms between
the lines are injective and the diagram commutes.
There are obvious different $T^2$-fibrations of
$T^4$. For $T^2 \times K^2$, diferent $T^2$ fibrations can be seen by
taking in 
$$\xymatrix{T^2 \times K^2 \ar[d] & \\ K^2 \times K^2 \ar[r] \ar[d] & K^2
\\ K^2 &}$$
the projections of $K^2 \times K^2$ onto different factors. 
\end{proof}

\begin{cor}\label{fib-cov-rel}
Any $\pi_1$-injective map $f: (E_1, \partial E_1) \to (E_2, \partial
E_2)$ of torus bundles over surfaces with non-empty 
$\pi_1$-injective boundary is homotopic rel boundary 
to a fiber-\-co\-ve\-ring map. 
\end{cor}

\begin{proof}
The condition on base implies that we are in the Case 1 of 
Proposition \ref{fiber-cover}. Hence $f$ is homotopic to a 
fiber-covering map. By the lemma \ref{q1}
it means that induced map on $\pi_1$'s send the fiber subgroup of $E_1$
into the fiber subgroup of $E_2$. From where 
$f \arrowvert_{\partial E_1}$ is homotopic to a fiber-covering map too,
because $\partial E_i$ is a subbundle of $E_i$.  

Denote the corresponding homotopies by 
$\{ f_t \} : E_1 \times I \to E_2$ and $\{ g_t \} : 
\partial E_1 \times I \to \partial E_2$. For each $t$, the step 
maps $f_t$ and $g_t$ are both homotopic to 
$f \arrowvert_{\partial E_1 \times \{ 0 \}}$. 
Presenting $E_1$ as $(\partial E_1 \times [0;t]) \cup 
(\partial E_1 \times [t;1]) \cup E_1$, one can define
a new homotopy $\{ H_t \} : E_1 \times I \to E_2$  of $f$ as follows.
On $\partial E_1 \times [0;t]$, the step map $H_t$ 
will be the (reparametrized) homotopy 
between $g_t$ and $f \arrowvert_{\partial E_1 \times \{ 0 \}}$ 
followed by the reparametrized homotopy between 
$f \arrowvert_{\partial E_1 \times \{ 0 \}}$ and $f_t$. 
On $(\partial E_1 \times [t;1]) \cup E_1 = E_1$, the map 
$H_t$ will be $f_t$.

The end map $H_1: (\partial E_1 \times [0;1]) \cup E_1 \to E_2$ 
is fiber covering on $E_1$ and on $(\partial E_1 \times \{ 0 \})$. 
Denote by $\gamma_0$ and $\gamma_1$ loops in $B_2$, subbundles over which
are covered by $H_1(\partial E_1 \times \{ 0 \})$ and  
$H_1(\partial E_1 \times \{ 1 \})$. As $\gamma_0$ and $\gamma_1$ are
homotopic in $B_2$, the homotopy between $H_1 \arrowvert_{\partial E_1
\times \{ 0 \}}$ and $H_1 \arrowvert_{\partial E_1 \times \{ 1 \}}$
can made fiber covering, and this new homotopy is homotopy to the old
one. Then $\{ H_t \}$ followed by the new homotopy gives a homotopy
{\it relatively to the boundary} between $f$ and a fiber covering map.
\end{proof}


\begin{cor}\label{he-rb}
Any homotopy equivalence rel boundary of torus bundles over surfaces 
with non-empty $\pi_1$-injective boundary is homotopic rel boundary 
to a diffeomorphism. 
\end{cor}
\begin{proof}
According to Corollary \ref{fib-cov-rel}, both maps 
of the homotopy equivalence can be made fiber covering maps 
by homotopy rel boundary. As both of them are of degree $\pm 1$, 
they are isomorphisms on the fibers. 
This and the commutative diagram of fundamental groups
imply that the monodromies are conjugate, hence the bundles are
isomorphic. As the obtained diffeomorphism of aspherical total spaces 
induces the same preserving peripheral structure isomorphism of
$\pi_1$'s as the initial map, they are homotopic rel boundary.
\end{proof}

Recall that a subgroup $A$ is said to be square root closed in $G$ if for
every element $g \in G$ such that $g^2 \in A$ one has $g \in A$,
too.

\begin{prop}\label{squareroot}
Let $B$ surface with boundary, $S$ is a component of $\partial B$, $E \to
B$ a $T^2$-bundle over $B$. Then the image $\pi_1(E \arrowvert_{S}) \to
\pi_1(E)$ is square root closed if and only if $B$ is not a M\"{o}bius band.
\end{prop}
\begin{proof}
If $B=D^2$, the homomorphism $\pi_1(E \arrowvert_{S}) \to \pi_1(E)$ is
onto and the statement is obvious. 

\noindent
If $B$ is different from the disc and M\"{o}bius band, $\pi_1(E
\arrowvert_{S}) \to \pi_1(E)$ is injective. As the diagram 
$$\xymatrix{ 0 \ar[r] &  \pi_1(E \arrowvert_{S}) \ar[r] \ar[d] &
\pi_1(E) \ar[d] \\
0 \ar[r] & \pi_1(S )  \ar[r] & \pi_1(B) }$$
\noindent commutes, $\pi_1(E \arrowvert_{S}) \subset \pi_1(E)$ is square root
closed if and only if $\pi_1(S) \subset \pi_1(B)$ does. Suppose 
$\pi_1(S) \subset \pi_1(B)$ is not square root closed.
Choose $a \notin \pi_1(S)$ with $a^2 \in \pi_1(S)$.
As $\pi_1(B)$ is free, and square roots are unique in free groups, so $a^2$
must be an odd power of the generator of $\pi_1(S)$.

Next observe that there is a  M\"{o}bius band $(M, \partial M) \to (B,S)$
with $\pi_1(M) \to a, \ \pi_1(\partial M) \to a^2$. Attach disks to $M$
and $B$ to get a map of $\mathbb{R}P^2=M \cup D^2 \to B \cup_{S} D^2$.
This induces $$\mathbb{Z}_2 \cong H_1(\mathbb{R}P^2; \mathbb{Z}_2) \to
H_2(B \cup_S D^2; \mathbb{Z}_2) \to H_2(D^2,S;\mathbb{Z}_2) \cong
\mathbb{Z}_2.$$
\noindent The composition is the same as boundary map $H_1(\partial M;
\mathbb{Z}_2) \to H_1(S; \mathbb{Z}_2)$ which is an isomorphism because
$a$ is an odd power of the generator. Therefore we conclude $H_2(B \cup_S
D^2; \mathbb{Z}_2 ) \cong \mathbb{Z}_2$ and $\mathbb{R}P^2 \to B \cup_S
D^2$ is an isomorphism on $H_2$ with $\mathbb{Z}_2$ coefficients. It
follows that $B \cup_S D^2$ is closed and $\pi_1(\mathbb{R}P^2) \to
\pi_1(B \cup_S D^2)$ has finite odd index. But $\mathbb{R}P^2$ is the
only closed surface with finite $\pi_1$, so
$B \cup_S D^2 \cong \mathbb{R}P^2$, and $B$ is a M\"{o}bius band.
\end{proof}

\section{Graph-manifolds}
We use the term
\emph{non-singular block} for the total space of a $T^2$-bundle over
a compact surface (with non-empty boundary) different
from a 2-disc, an annulus and a M\"{o}bius band (hence, a
surface with free non-abelian fundamental
group). Boundary components of blocks are $T^2$-bundles over
$S^1$ and are $\pi_1$-injective in blocks.

\begin{defi}\label{defigraph}\rm
A 4-dimensional closed connected compact oriented manifold
is a \textit{non-singular graph-manifold} if it can be obtained by gluing
several blocks by diffeomorphisms of their boundaries.
\end{defi}

\noindent Fot simplicity we will say ``blocks'' instead of 
``non-singular blocks''
and ``graph-manifolds'' instead of ``non-singular graph-manifolds''.

\noindent {\bf Example}\qua The simplest examples of $4$-dimensional
graph-manifolds are $T^2$-bundles over closed hyperbolic surfaces
(all the glueing diffeomorphisms being
trivial). A more interesting examples can be
constructed by taking oriented $S^1$-bundles over some
$3$-dimensional graph-manifolds: for instance, such that all their
blocks have $\pi_1$-injective boundary components (for exemple,
lens spaces are not good as bases) and all the blocks being
locally trivial $S^1$-bundles (i.e.\ no exceptional Seifert
fibers).

Any decomposition as a union of blocks will be called a
\emph{graph-structure}.
Topologically, a graph-structure is determined by a
system of embedded $\pi_1$-injective $T^2$-bundles over circles,
called \emph{decomposing ma\-ni\-folds}.
A graph-structure is
\emph{reduced} if all the glueing diffeomorphisms are not
fiber-\-preserving, or, equivalently, if the induced isomorphism
of $\pi_1$'s does not preserve the fiber subgroup.
As the fiber subgroup is unique
in $\pi_1$ of the block, the notion of reduced structure is well
defined.

{\bf Immediate properties}\qua Graph-manifolds are
aspherical: since inclusions of boundary components into blocks are
$\pi_1$-injective, they are graphs of aspherical spaces,
and the universal covering
of a graph of aspherical spaces is contractible (\cite{SW},
prop.~3.6 p.156). The definition also implies that the Euler
characteristic of graph-manifolds is $0$, because the Euler 
characteristic of block is $0$, and gluings are made along
3-manifolds. Finally, graph-manifolds can be smoothed:
the blocks are smooth, and gluings are made by diffeomorphisms of 3-manifolds.
More, a given graph-structure determines the smoothing in a unique
way, because the smooth structure on a $3$-manifold is unique
and homotopic diffeomorphisms of torus bundles over $S^1$ are
isotopic \cite{Wald}.

\begin{prop}
The signature of a closed oriented graph-manifold $W^4$ with
reduced graph-structure all the blocks of which have orientable
bases is $\ \sigma(W^4)=0 .$
\end{prop}
\begin{proof}
The blocks of an orientable graph manifold are orientable, and the
signature of the graph-manifold induces the orientations on the blocks. One
can assume that all the orientations of blocks are such that glueing
diffeomorphisms reverse the induced orientation of boundaries.
The orientation on a block comes from the orientation of
its fiber plus the orientation of its base. Hence we can speak about
presentation of boundaries of blocks as some 
$M_{\varphi_i}=(T^2 \times I) / (x;0) \sim (\varphi_i(x); 1), 
\ \varphi_i \in SL(2, \mathbb{Z})$.

Determine first the signatures of blocks. In a reduced graph-structure,
the boundaries of all the blocks have many non-isotopic
$T^2$-bundle structures. Hence the monodromies of all decomposing
manifolds must be conjugate to $\bigl(
\begin{smallmatrix} 1 & n_i\\ 0 & 1
\end{smallmatrix} \bigr) $ (Proposition \ref{q2}). Thus by Meyer's Theorem \cite{Meyer},
the signature of a block $M^4$ of such a manifold is {\small
$$\sigma (M^4) = \frac{1}{3} \sum_{i=1}^{k} n_i .$$}Now 
apply the Novikov's additivity to get $\sigma(W^4)$ 
by adding the signatures of the
blocks. When one glues the blocks, $M_{\varphi}$ can be glued either with $M_{a
\varphi a^{-1}}$ or with $M_{a \varphi^{-1} a^{-1}}$ for some $a \in SL(2,
\mathbb{Z})$. As Meyer's characteristic function $\Psi: SL(2, \mathbb{Z})
\to \mathbb{Z}$ is invariant under conjugation in $SL(2, \mathbb{Z})$, for
the signature calculation one can assume that $M_{\varphi}$ can be glued
either with $M_{\varphi}$ or with $M_{\varphi^{-1}}$.
There exists an orientation reversing diffeomorphism $M_{\varphi} \to
M_{\varphi}$ if and only if the Euler number of $S^1$-bundle of
$M_{\varphi}$ with a triangular $\varphi$ is $0$ \cite{Seifert}, i.e.\ if
$\varphi = \bigl( \begin{smallmatrix} \pm 1 & 0\\ 0 & \pm 1
\end{smallmatrix} \bigr) $. But the boundary component with such
monodromy gives contribution $0$ into the signature.
There always exists an orientation reversing diffeomorphisms $M_{\varphi}
\to M_{\varphi^{-1}}$, but $\Psi \bigl( \begin{smallmatrix}
 1 & n\\ 0 & 1 \end{smallmatrix} \bigr) +
\Psi \bigl( \begin{smallmatrix}
 1 & -n \\ 0 &  1 \end{smallmatrix} \bigr)=n + (-n)=0$ and
$\Psi \bigl( \begin{smallmatrix}
- 1 & \ n\\ 0 & - 1 \end{smallmatrix} \bigr) +
\Psi \bigl( \begin{smallmatrix}
- 1 & -n \\ 0 & -1 \end{smallmatrix} \bigr) =n + (-n)=0$, hence every such
pair also gives the contribution $0$ in the signature.
Hence after adding the signatures of all the blocks we will
obtain $\sigma (W^4)=0$.
\end{proof}

\begin{lemma}\label{into-bound}
Let $M$ be a $T^2$-bundle over surface with non-empty $\pi_1$-injective
boundary. Then any non-\-fi\-ber-\-co\-ve\-ring $\pi_1$-injective map
$$f: (T^2 \times I; \partial(T^2 \times I)) \to (M; \partial M)$$
\noindent sending $\partial(T^2 \times I)$ into the same boundary
component is homotopic rel boundary to a map into $\partial M$.
\end{lemma}
\begin{proof}
Denote the component of $\partial M$ containing
$f( T^2 \times \partial I)$ by $M_{\varphi}$.
Fix a point $t \in T^2$. As $T^2 \times I$ is aspherical, we have to show that
$f \arrowvert_{t \times I}: (I,\partial I) \to (M, M_{\varphi})$
is homotopic to a map into $M_{\varphi}$.

Denote the projections $p: M \to B^2$, $p_b: M_{\varphi} \to S^1$
and natural inclusions $l: M_{\varphi} \rightarrow M$,
$l_p: p(M_{\varphi}) \rightarrow B^2$. Take a path
in $M_{\varphi}$ that joins $f(t \times \{ 0 \})$ and
$f(t \times \{ 1 \})$; the union of this path with $f (t \times I)$
is an element of $\pi_1(M;b)$ where $b=f(t \times \{ 0 \})$.

In the diagram
$$\xymatrix{ 0 \ar[rd]&&&&
\\
&  \mathbb{Z} \oplus \mathbb{Z} \ar[rd]^{f_{\ast}} & 0 \ar[d]
&  0 \ar[d] &
\\
0 \ar[r] & \mathbb{Z} \oplus \mathbb{Z} \ar[r]
\ar@{=}[d] & \pi_1(M_{\varphi}) \ar[r]^{{p_b}_{\ast}}
\ar[d]^{l_{\ast}}  & \mathbb{Z} \ar[r]
\ar[d]^{{l_p}_{\ast}} & 0 \\
 0 \ar[r] & \mathbb{Z} \oplus \mathbb{Z} \ar[r] &
\pi_1(M) \ar[r]^{p_{\ast}} & \pi_1(B^2) \ar[r] & 0 }$$
denote $G=Im \, f_{\ast}$; ${p_b}_{\ast}(G)$ is non-trivial (by
assumption). We have $\gamma l_{\ast} f_{\ast}(G) \gamma^{-1}$ $\subset Im
\, l_{\ast},$ hence $$p_{\ast}(\gamma) \, {l_p}_{\ast} {p_b}_{\ast}(G) \,
p_{\ast}(\gamma^{-1}) \subset Im \, ({l_p}_{\ast} {p_b}_{\ast}) \cong
\mathbb{Z}.$$
\noindent As $ {l_p}_{\ast} {p_b}_{\ast}(G) $ is abelian and
non-trivial, and $ Im \, ({l_p}_{\ast} {p_b}_{\ast}) $ is generated by
primitive element of $\pi_1(B^2)$, we conclude that $p_{\ast}(\gamma)
\in Im \, ({l_p}_{\ast})$ \cite{LyndSch}, hence $\gamma \in Im \, l_{\ast}$.
\end{proof}

\begin{prop}\label{main-prop}
Let $W=\cup W_i$ and $W'=\cup W'_k$ be non-singular graph-man\-if\-olds with
reduced graph-structures. Then any $\pi_1$-injective map $f:W \to W'$ is
homotopic to $\bigcup f_i$, where each $f_i : (W_i, \partial W_i) \to (W'_j,
\partial W'_j)$ is fiber covering map.
\end{prop}
\begin{proof}$\phantom{99}$

{\bf Step~1}\qua {\sl Any $\pi_1$-injective map of torus bundle over circle
$f: M_{\varphi} \to W'$ is homotopic to a fiber covering map into one block.}

By Main Technical Result, one can move $f$ by homotopy such that 
the inverse image by $f$ of decomposing submanifolds
becomes disjoint union of $\pi_1$-injective $2$-tori, 
embedded in $M_{\varphi}$. Then
$M_{\varphi}$ cut along them is either $\bigcup_{i} (T^2 \times I)_i$ or
$\bigcup_{i} (T^2 \times I)_i \cup (K^2 \widetilde{\times} I) 
\cup (K^2 \widetilde{\times} I) $, each summand lying in a block 
($K^2 \widetilde{\times} I$ is twisted oriented $I$-bundle over Klein bottle).

{\bf Case 1}\qua $M_{\varphi}$ cut along preimages of decomposing
submanifolds is $\bigcup_i (T^2 \times I)$.

\noindent 
First observe that if $f$ sends $T^2 \times \{ 0 \}$ and $T^2 \times \{ 1 \}$
into different boundary components of a block,
then $f$ is homotopic to a fiber-covering map.
Indeed, denote the block $p: M \to B$ and choose a base point $b
\in f(T^2 \times \{ 0 \})$. Then $p_{\ast} f_{\ast} (\mathbb{Z} \oplus
\mathbb{Z})$ lies in the subgroup of $\pi_1(B, p(b))$ corresponding to
the component of $\partial B$ containing $pf(T^2 \times \{ 0 \})$, and
in the conjugation class of the subgroup corresponding to
the component of $\partial B$ containing $pf(T^2 \times \{ 1 \})$. 
But conjugate classes of different boundary components can intersect only
if $B$ is an annulus, because the conjugation defines a map $S^1 \times I
\to B$ of non-zero degree. Hence $p_{\ast} f_{\ast}
(\mathbb{Z} \oplus \mathbb{Z})=1$ which means that $f$ sends the fiber 
subgroup of $T^2 \times I$ into the fiber subgroup of the block.
Hence $f$ is homotopic to fiber-covering map. 

\medskip
{\bf Remark~1}\qua This observation implies that
in a reduced graph-structure the fibers of different blocks are not
homotopic. Indeed, if a graph-manifold has just 2 blocks, then the claim
comes from the definition of the reduced graph-structure.
If there are more blocks, take one of them, its fiber satisfies 
the conditions of the previous observation in all
the neighboring blocks. Hence, in every neighboring block this fiber
is not homotopic to any torus in the remaining
boundary components. But in these components lie in
particular the fibers of the next neighboring blocks etc.
\medskip

\noindent As the graph-structure is
reduced, in all the neighboring blocks $f$ is not homotopic to a
fiber-covering map and hence $f(\partial(T^2 \times I))$ lie in the same
boundary component. Hence one can apply Lemma
\ref{into-bound} to the neighboring $(T^2 \times I)$'s and move them into
the block where the fiber-covering $f(T^2 \times I)$ lies.

\begin{figure}[ht]\label{}
\centering
\input{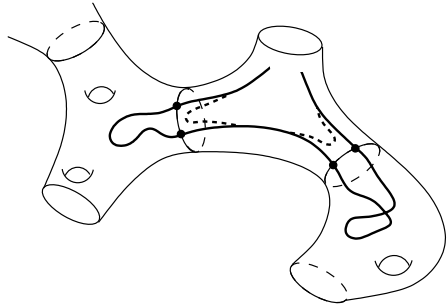}
\protect\caption{{\small Shrinking in one block}}
\end{figure}

{\bf Case 2}\qua $M_{\varphi}$ cut along preimages of decomposing
submanifolds is $\bigcup_i (T^2 \times I) \cup (K^2 \widetilde{\times} I)
\cup (K^2 \widetilde{\times} I) $.

\noindent Take the first copy of 
$K^2 \widetilde{\times} I, \ I=[ -1; 1]$, denote the
decomposing manifold, in which the image under $f$ of its boundary
lies, by $M_{\varphi_1}$. Its boundary torus is a two-fold
covering of the Klein bottle in the base and if $\pi_1( K^2
\widetilde{\times} I) = \pi_1( K^2)=\langle a, b \arrowvert
aba^{-1}=b^{-1} \rangle$, then the boundary torus corresponds to
the subgroup $  \langle a, b^2 \rangle$.
As the subgroups of boundary components are square root closed 
in the fondamental groups of the blocks (Proposition \ref{squareroot}), 
the subgroup $f_{\ast} (\pi_1 ( K^2 \widetilde{\times} I), x)$ must lie 
in the subgroup corresponding to $M_{\varphi_1}$, because\break
$f_{\ast} (\pi_1 (\partial ( K^2 \widetilde{\times} I)), x)$ lies there.
As $K^2 \widetilde{\times} I$ is aspherical, 
$f \arrowvert_{K^2 \widetilde{\times} I}$ can be moved by
homotopy in $M_{\varphi_1}$ and, hence, out of its original block
in the neighboring one. Repeat the previous reasonnings for the
union of $(K^2 \widetilde{\times} I)$ with the next $T^2 \times I$
gives a new $(K^2 \widetilde{\times} I)$.
In the end it will be two copies of $(K^2 \widetilde{\times} I)$, 
and the image of each of them
under $f$ can be moved into the same decomposing manifold. Hence, in 
this case $f$ is homotopic to a map into a
decomposing manifold.
 
Once $f(M_{\varphi})$ is shrinked in one block,
look at the homomorphism that $f$ induces on $\pi_1$'s. The image of the 
subgroup of the fiber of $M_{\varphi}$ vanishes when projecting on 
$\pi_1$ of the base of $M$, because it is abelian normal subgroup of 
non-abelian free group. Hence, by lemma \ref{q1} the map is homotopic 
to a fiber covering one.

{\bf Step~2}\qua {\sl Any $\pi_1$-injective map of a block $f: M \to W'$ is
homotopic to a map into one block of $W'$.}

Any block retracts on a torus bundle over a wedge of circles. Torus
bundle over a wedge of circles can be obtained from a torus bundle over
circle (with the monodromy equal to the product of the monodromies of
the petals) by identifying some fibers. Change the map of this ``big'' single
torus-bundle by homotopy given by Step~1. The images of fibers that 
are identified are two by two homotopic by two kinds of homotopy.
The first are given by petals and now lies in one block $M'$. The second comes
from identification and still lies in the whole $W'$. In order to
identify the images, one has to shrink these homotopies into $M'$.

Apply the Main Technical Result to each of these homotopies, this will make
their intersections with decomposing submanifolds $W'$
$\pi_1$-injective.

If the tori that must be identified are fiber-covering in $M'$, the
corresponding homotopies can lie only in the neighboring blocks,
where they are not-fiber-covering, hence by Lemma \ref{into-bound} can be
shrinked in $M'$. 

If the tori that must be identified are not 
fiber-covering in $M'$, then the homotopies lie in the union
of $M'$ with its neighboring blocks, because if the part of the
homotopy in the neighboring block does fiber-covering, then in the
following blocks it does not and Lemma \ref{into-bound} does apply.
More, all this homotopies must lie in just one neighboring block of $M'$, 
because elsewhere we would have a $\pi_1$-injective map 
$(T^2 \times I, \partial (T^2 \times I)) \to (M', \partial M'))$ which
would be non-fiber-covering but sending $\partial (T^2 \times I)$ into
different components of $\partial M'$. Then the remaining parts of
homotopies and the images of all the petals can be shrinked in this
neighboring block.  

{\bf Step~3}\qua For all $i$, change $f \arrowvert_{W_i}$ according 
to Step~2, denote it by $f_i$. Let $W_i, W_k \subset W$ be 
neighboring blocks. For each component of $W_i \cap W_k$, by Step~1, the
block in which it lies is unique. 
As $f_i \arrowvert_{W_i \cap W_k}$ and $f_k \arrowvert_{W_i \cap W_k}$ 
are homotopic, one conclude that the corresponding $W'_i$ and $W'_k$ are
neighboring and $f_i \arrowvert_{W_i \cap W_k}, 
f_k \arrowvert_{W_i \cap W_k}$ are homotopic to a map into 
$W'_i \cap W'_k$. Use the corresponding homotopy inside $W'_i$ 
to define the new $f_i : W_i=W_i \bigcup ((W_i \cap W_k) \times I) \to W'_i$
being the homotopy on $(W_i \cap W_k) \times I$ part and the old 
$f_i$ on $W_i$ part; then make the same for $f_k$. 
By doing it on all pairs of neighboring blocks, one obtain the map
$f= \bigcup f_i$ with $f_i : (W_i, \partial W_i) \to (W'_i,
\partial W'_i)$. Apply Corollary \ref{fib-cov-rel} to every $f_i$ to make
it fiber covering. It remains to bind the new 
$f_i \arrowvert_{W_i \cap W_k}$ and 
$f_k \arrowvert_{W_i \cap W_k}$ by a fiber covering homotopy inside the 
correponding $W'_i \cap W'_k$, which is possible because they covers the
same subbundles of $W'_i \cap W'_k$ with the same degree. 
We obtain a map $f = \bigcup f_i: W \to W'$ such that every 
$f_i : (W_i, \partial W_i) \to (W'_i, \partial W'_i)$ is fiber covering.
\end{proof}

\begin{theorema}
Any homotopy equivalence between non-singular graph-manifolds with reduced
graph-structures is homotopic to a diffeomorphism.
\end{theorema}
\begin{proof}
Take the collars of decomposing manifolds in the blocks, the blocks
without this collars remain the blocks. For each
decomposing manifold $M_{\varphi}$, the union of its collars on both side 
is $M_{\varphi} \times I$, call it
``double collar of $M_{\varphi}$'' in the graph-manifold.

According to Proposition \ref{main-prop}, both maps of the homotopy 
equivalence $f: W \to W', g: W' \to W$ can be moved by homotopy such that
their restrictions on blocks without collars are 
rel boundary fiber covering maps. 
As after homotopies we still have $gf \sim id_{W}$, 
the fiber of a block (without collars) $M$ is homotopic to its image
by $gf$. Hence, by Remark~1, $gf(M) \subset M$ for every block of $W$,
i.e.\ the restrictions of $f$ and $g$ gives the homotopy equivalences rel
boundary of blocks without collars. 
According to Corollary \ref{he-rb}, these resrictions
are homotopic rel boundary to diffeomorphisms. One has to bind the
obtained block's diffeomorphisms on the double collars of the decomposing 
manifolds. For this note that the diffeomorphisms, that are
the restricions of $f$ and $g$ on the boundaries of blocks without
collars, are homotopic. Hence restrictions of $f$ and $g$ on double
collars of decomposing manifolds are homotopy equivalences rel boundary
that are diffeomorphism on the boundaries. As decomposing manifolds are
sufficiently large, these homotopy equivalences are homotopic rel boundary
to diffeomorphisms, by homotopies that are constant on the boundaries.
\cite{Wald}.
\end{proof}

{\bf Acknowledgements}\qua
This work has been done as a thesis in {\it co-tutelle} under supervision of 
Claude Hayat-Legrand and Vladimir Sharko, to whom the author is very grateful. 
The author also wishes to express gratitude to Alexis Marin for 
fruitful criticism. It is also a pleasure to
thank Frank Quinn for very helpful discussions.

\Addresses\recd

\end{document}